\documentclass[12pt]{article}

\usepackage{graphicx,color,ulem,dsfont}
\usepackage{hyperref}
\usepackage{epsfig, tkz-berge,tikz,pgf}
\usepackage{amsmath}

\newtheorem{thm}{Theorem}[section]
\newtheorem{prop}[thm]{Proposition}
\newtheorem{cor}[thm]{Corollary}
\newtheorem{lem}[thm]{Lemma}
\newtheorem{conj}[thm]{Conjecture}
\newtheorem{exa}[thm]{Example}
\newtheorem{defn}[thm]{Definition}
\newtheorem{rem}[thm]{Remark}
\newtheorem{que}[thm]{Question}
\newtheorem{obs}[thm]{Observation}
\newtheorem{prob}[thm]{Problem}

\newcommand{\ben}{\begin{enumerate}}
\newcommand{\een}{\end{enumerate}}
\newcommand{\ble}{\begin{lem}}
\newcommand{\ele}{\end{lem}}
\newcommand{\bth}{\begin{thm}}
\newcommand{\eth}{\end{thm}}
\newcommand{\bpr}{\begin{prop}}
\newcommand{\epr}{\end{prop}}
\newcommand{\bco}{\begin{cor}}
\newcommand{\eco}{\end{cor}}
\newcommand{\bcon}{\begin{conj}}
\newcommand{\econ}{\end{conj}}
\newcommand{\bde}{\begin{defn}}
\newcommand{\ede}{\end{defn}}
\newcommand{\brem}{\begin{rem}}
\newcommand{\erem}{\end{rem}}
\newcommand{\bque}{\begin{que}}
\newcommand{\eque}{\end{que}}
\newcommand{\bobs}{\begin{obs}}
\newcommand{\eobs}{\end{obs}}
\newcommand{\bprob}{\begin{prob}}
\newcommand{\eprob}{\end{prob}}
\newcommand{\bex}{\begin{exa}}
\newcommand{\eex}{\end{exa}}
\newcommand{\barr}{\begin{array}}
\newcommand{\earr}{\end{array}}
\newcommand{\btab}{\begin{tabular}}
\newcommand{\etab}{\end{tabular}}
\newcommand{\beq}{\begin{equation}}
\newcommand{\eeq}{\end{equation}}
\newcommand{\bea}{\begin{eqnarray*}}
\newcommand{\eea}{\end{eqnarray*}}
\newcommand{\bce}{\begin{center}}
\newcommand{\ece}{\end{center}}
\newcommand{\bpi}{\begin{picture}}
\newcommand{\epi}{\end{picture}}
\newcommand{\bfi}{\begin{figure} \begin{center}}
\newcommand{\efi}{\end{center} \end{figure}}

\newcommand{\bsl}{\begin{slide}{}}

\newcommand{\pf}{{\bf Proof.}}
\newcommand{\qed}{\rule{1ex}{1ex}}

\begin{document}

\title{A Graph Theoretical Analysis of the Number of Edges in $k$-Dense Graphs}

\author{
 {\bf Linda Eroh$^1$},
 {\bf Henry Escuadro$^2$},
  \and
 {\bf Ralucca Gera$^3$},
 {\bf Samuel Prahlow$^4$},
 and
 {\bf Karl Schmitt$^5$}\\ \\
 $^1$\small Department of Mathematics, University of Wisconsin-Oskosh,\\
\small Oshkosh, WI $54901$; {\small \em eroh@uwosh.edu}\\
$^2$\small Department of  Mathematics, Juniata College,\\
\small Huntingdon, PA $16652$; {\small \em escuadro@juniata.edu}\\
$^3$\small Department of Applied Mathematics,  Naval Postgraduate School, \\
 \small
 Monterey, CA $93943$;  {\small \em rgera@nps.edu}\\
$^4$\small Mathematics and Statistics Department,  Valparaiso University, \\
 \small
Valparaiso, IN $46383$;  {\small \em samuel.prahlow@valpo.edu}\\
$^5$\small Mathematics and Statistics Department,  Valparaiso University, \\
 \small
Valparaiso, IN $46383$;  {\small \em karl.schmitt@valpo.edu}\\
}

\maketitle

\begin{abstract}
Due to the increasing discovery and implementation of networks within all disciplines of life, the study of subgraph connectivity has become increasingly important. Motivated by the idea of community (or sub-graph) detection within a network/graph, we focused on finding characterizations of k-dense communities. For each edge $uv\in E(G)$, the {\bf edge multiplicity} of $uv$ in $G$ is given by $m_G(uv)=|N_{G}(u)\cap N_{G}(v)|.$ For an integer $k$ with $k\ge 2$, a {\bf $k$-dense community} of a graph $G$, denoted by $DC_k(G)$, is a maximal connected subgraph of $G$ induced by the vertex set
$$V_{DC_k(G)} = \{v\in V(G) : \exists  u\in V(G)\ {\rm such\ that\ } uv\in E(G)\ {\rm and\ } m_{DC_{k(G)}}(uv)\ge k-2\}.$$
In this research, we characterize which graphs are $k$-dense but not $(k+1)$-dense for some values of $k$ and study the minimum and maximum number of edges such graphs can have. A better understanding of $k$-dense sub-graphs (or communities) helps in the study of the connectivity of large complex graphs (or networks) in the real world.

\end{abstract}

\vspace{.1cm}

\noindent
{\bf Key Words:}  $k$-dense subnetworks (or $k$-dense subgraph), $k$-dense community, $k$-dense graph, $k$-core, $k$-core subnetwork

\vspace{.2cm}

\noindent
{\bf AMS Subject Classification:} 05C69, 05C75, 05C62.


\section{Motivation and definitions}  \label{section:motivation}

This section covers definitions and notation that will be relevant for the whole study. Please refer to \cite{clz} for notation and terminology that are not discussed in this paper. In this study, all graphs are simple (i.e., no loops and no multiple edges). Recall that $\delta(G)$ denotes the minimum degree and $\Delta(G)$ denotes the maximum degree among all vertices in $G$.  Also, a graph $G$ is $k$-edge connected, if for any set $S$ of fewer than $k$ edges $G-S$ is connected, and there is a set $T$ of $k$ edges such that $G-T$ is disconnected.

\medskip

Also recall that a clique in a graph $G$ is a complete subgraph of $G$.  The order of a largest clique is the clique number $\omega(G).$  Note that the maximum clique in a graph may not be unique.  It was shown by  Karp in \cite{K} that the decision problem of finding a clique of some integer size $k$ in a given graph $G$ is NP-complete.

\medskip

For a graph $G$ and vertex $v \in V(G)$, the open neighborhood of $v$, denoted by $N_G(v)$, is the set of all neighbors of $v$, namely,
$N_G(v) = \{u : uv \in E(G)\}$.  If we restrict our attention to a subgraph $H$ of $G$,  we then have $N_H(v) = \{u : uv \in E(H)\}$.  The closed neighborhood of $v$ is
$N[v] = N(v) \cup \{v\}$.

\medskip

In this paper, we study the idea of $k$-dense subgraphs, introduced by Saito, Yamada and Kazama in \cite{SYK} as an alternative to cliques in detecting communities in a graph.  We believe this is an alternative to possibly detecting cliques in graphs, since the inferred graphs may not always contain the complete information, such inferring all router connections in a network.  Thus, a clique in the real topology may be a $k$-dense community in the inferred and studied graph.  The following definitions were introduced in \cite{SYK} but we rephrase them here in graph theoretical terms.

\begin{defn}
\label{multiplicity_kdense}
Let $G$ be a graph. For each edge $uv\in E(G)$, the {\bf edge multiplicity} of $uv$ in $G$ is given by $m_G(uv)=|N_{G}(u)\cap N_{G}(v)|.$
\end{defn}


\begin{defn}
For an integer $k$ with $k\ge 2$, a {\bf $k$-dense community} of a graph $G$, denoted by $DC_k(G)$, is a maximal connected subgraph of $G$ induced by the vertex set
$$V_{DC_k(G)} = \{v\in V(G) : \exists  u\in V(G)\ {\rm such\ that\ } uv\in E(G)\ {\rm and\ } m_{DC_{k(G)}}(uv)\ge k-2\}.$$
\end{defn}

In other words, $DC_k(G)$ is a connected subgraph of $G$ with $k$ or more vertices in which every two adjacent vertices have at least $k-2$ neighbors in common in $DC_{k}(G)$.  Note that a $k$-dense community might contain a connected proper subgraph $G[S]$ where  $S\subseteq V(G)$  and for every pair of vertices $u,v$ in $S$ that are adjacent in $G[S]$, we have $m_{G[S]}(uv)\ge k-2.$ Here, $G[S]$ denotes the subgraph of $G$ induced by $S.$  In this case, we call $S$ a {\bf $k$-dense sub-community} of $G.$     For example, consider the graph $G$ obtained from two copies of $K_5$, by identifying 3 of the vertices of both copies (thus $G=K_3 + 2K_2$).  Each of the original copies of $K_5$ is a $5$-dense sub-community of $G$ while the graph $G$ itself is a $5$-dense community. 

\medskip

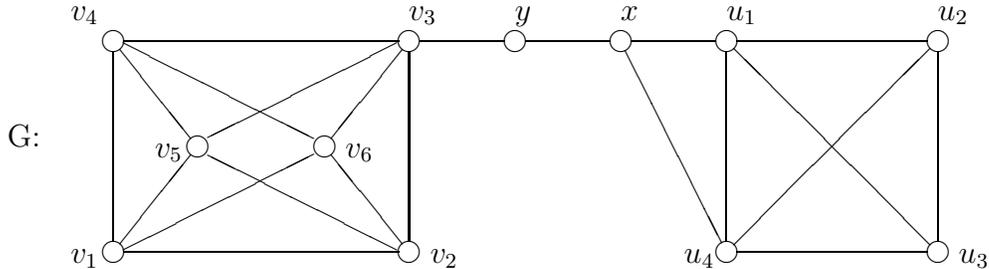
\begin{figure}[h]~\label{Fig:dense}
\setlength{\unitlength}{4pt}
\begin{center}
\begin{picture}(80,20)(0,0)
\put(-10, 10){G:}
\put(0,0){\circle{2}}
\put(-4,-1){$v_1$}
\put(8,10){\circle{2}}
\put(4,9){$v_5$}
\put(20,10){\circle{2}}
\put(22,9){$v_6$}
\put(28,0){\circle{2}}
\put(30,-1){$v_2$}
\put(0,20){\circle{2}}
\put(-4,22){$v_4$}
\put(28,20){\circle{2}}
\put(28,22){$v_3$}
\put(1,0){\line(1,0){26}}
\put(1,20){\line(1,0){26}}
\put(0,1){\line(0,1){18}}
\put(28,1){\line(0,1){18}}
\put(.63,.8){\line(4,5){6.74}}
\put(.63,19.2){\line(4,-5){6.74}}
\put(27.37,.8){\line(-4,5){6.74}}
\put(27.37,19.2){\line(-4,-5){6.74}}
\put(.97,.24){\line(2,1){18}}
\put(.97,19.76){\line(2,-1){18}}
\put(26.84,.24){\line(-2,1){17.9}}
\put(26.84,19.76){\line(-2,-1){17.9}}
\put(38,20){\circle{2}}
\put(38,22){$y$}
\put(48,20){\circle{2}}
\put(48,22){$x$}
\put(58,20){\circle{2}}
\put(58,22){$u_1$}
\put(78,20){\circle{2}}
\put(78,22){$u_2$}
\put(58,0){\circle{2}}
\put(54,-1){$u_4$}
\put(78,0){\circle{2}}
\put(80,-1){$u_3$}
\put(29,20){\line(1,0){8}}
\put(39,20){\line(1,0){8}}
\put(49,20){\line(1,0){8}}
\put(59,20){\line(1,0){18}}
\put(59,0){\line(1,0){18}}
\put(48.44,19.11){\line(1,-2){9.12}}
\put(58,1){\line(0,1){18}}
\put(78,1){\line(0,1){18}}
\put(58.71,.71){\line(1,1){18.58}}
\put(77.29,.71){\line(-1,1){18.58}}
\end{picture}
\caption{$G[\{v_i,u_j|1\leq i \leq 6, 1 \le j \le 4\}]$ is the $4$-dense subgraph of $G$, $G[\{v_j|1\leq j \leq 4\}]$ is a $4$-dense community that does not contain a $4$-clique, $G[V(G)-\{y\}]$ is the $3$-dense subgraph of $G$, while $G$ itself is $2$-dense. }
\end{center}
\end{figure}

As was shown in \cite{SYK},  a $k$-clique in a graph $G$ is also a $k$-dense sub-community in $G$.  However, the converse is not true. Consider the graph $H=G[\{v_1,v_2,\ldots, v_6\}]=C_4+2K_2$ shown in the Figure~\ref{Fig:dense}. Note that $H$ is a 4-dense community in $G$ that is not $K_4$ nor contains $K_4$ as a subgraph.

\medskip

One can observe that just like cliques in graphs, a given graph might have more than one $k$-dense community.  In this paper, we say that the union of all $k$-dense communities of $G$ is the {\it $k$-dense subgraph of $G$}, a concept that was introduced in  \cite{SYK} under the name of $k$-dense subnetwork.


\begin{defn}\label{DEFN:kdense}
For an integer $k\ge 2$, the {\bf $k$-dense subgraph} of a graph $G$, denoted by $D_k(G)$, is the union of all the $k$-dense communities in $G$. If $D_k(G)=G$, then we say that $G$ is a {\bf $k$-dense graph} or simply that $G$ is {\bf $k$-dense}. If $D_k(G)=G$ and $D_{k+1}(G)\ne G$, then we say that $G$ is a {\bf $k^*$-dense graph} or simply that $G$ is {\bf $k^*$-dense}.
\end{defn}

Note that the $k$-dense subgraph, $D_k(G)$, of a graph $G$ is unique but need not be connected. Also, a graph $G$ is $k^*$-dense if and only if $G$ is $k$-dense but not $(k+1)$-dense.

\medskip

In \cite{SYK}, Saito et al. introduced the idea of the $k$-dense subgraph of a graph $G$ as an alternative way to detect communities in graphs with many vertices and edges (viewed as large-scale complex networks). Earlier methods used $k$-cliques (complete subgraphs that contain $k$ vertices) and $k$-cores (maximal induced subgraphs whose minimum degree is $k-1$). By applying their algorithm to detect close-knit communities using the $k$-dense method on various real-life networks, Saito et al. found that their method is almost as efficient as the $k$-core method. Moreover, the communities found using the $k$-dense method proved to be comparable to those found using the $k$-clique method. Saito et al. indicate that a hierarchy of $k$-clique $\subseteq$ $k$-dense $\subseteq$ $k$-core exists, without providing a formal proof. They also state that the complete graph on $k$ vertices is $k$-dense (and a $k$-core) and that if $G$ is a $k$-dense graph with $k$ vertices, then $G$ must be a $k$-clique (complete). We restate the last claim as Proposition~\ref{PROP:ndense} and present a proof for it.

\medskip

While the authors of \cite{SYK} analyzed case studies on the concept of the $k$-dense subgraph, this paper pursues a graph theoretical study on the topic. In particular, we focus on determining lower and upper bounds for the number of edges in a graph with $n$ vertices that is $k^*$-dense. In Section~\ref{section:general} we present general bounds and a realization result on the number of edges in $k^*$-dense graphs.  In Section~\ref{section:k=3,k=4} we present the  minimum number of edges in $2^*$-dense, $3^*$-dense and $4^*$-dense graphs,
followed by the  maximum number of edges in a $k^*$-dense graph in Section~\ref{section:upper_bound}.  Finally, in Section~\ref{Section:complete k=3, 4}, we present complete results for $2^*, 3^*$-dense and $4^*$-dense graphs.


\section{General observations}  \label{section:observations}

Using the definitions from Section~\ref{section:motivation}, we have the following quick proposition.

\bpr\quad A graph $G$ is $k$-dense if and only if every edge of $G$ appears in at least $k-2$ different triangles. \epr

The next result tells us that the minimum degree of a $k$-dense graph is at least $k-1.$

\bpr \label{Prop:2.2}
\quad Let $G$ be a nontrivial graph that is $k$-dense. If $v \in V(G)$, then $\deg(v) \geq k-1$.  Furthermore, if $\deg(v) = k-1$, then $v$ is in a subgraph that is isomorphic to $K_k$.\epr

\noindent
\pf\quad
Suppose $v \in V(G)$ for some nontrivial graph $G$ that is $k$-dense.  Since $G$ is nontrivial and has no isolated vertices, $v$ must be incident with an edge $uv$.  Since $u$ and $v$ have at least $k-2$ common neighbors, it follows that $\deg(v) \geq k-1$.  Now, suppose $\deg(v) = k-1$, say $N(v) = \{v_1, v_2, \ldots, v_{k-1}\}$.  Then for each $i$, $1 \leq i \leq k-1$, the vertices $v$ and $v_i$ have $k-2$ common neighbors.  Since $\deg(v) = k-1$, those common neighbors must be $v_j$ for $1 \leq j \leq k-1$, $j \neq i$.  Thus, $v, v_1, v_2, \ldots, v_{k-1}$ induces a complete subgraph. \hfill\qed

\medskip

If $G$ is a $k$-dense graph, then $G-v$ is not necessarily $k$-dense. The following result tells us that $G-v$ is at least $(k-1)$-dense.

\bpr\quad Suppose $G$ is a $k$-dense graph, where $k \geq 3$, and $v \in V(G)$.  Then $G - v$ is at least $(k-1)$-dense. \epr

\noindent
\pf\quad Suppose $G$ is $k$-dense and $v \in V(G)$.  Let $u, w \in V(G-v)$ such that $uw \in E(G-v)$.  Then $u$ and $w$ are adjacent in $G$ as well. Hence, $u$ and $w$ have $k-2$ common neighbors in $G$,  one of which may be $v.$ This means that $u$ and $w$ will have at least $k-3$ common neighbors in $G-v$.  \hfill\qed

\medskip

 To see that the graph $G-v$ need not be $(k-1)^*$-dense if $G$ is $k^*$-dense, consider the graph $G$ obtained from a $5$-clique, whose vertices are $v_1, v_2, \ldots, v_5$, and a $6$-clique, whose vertices are $u_1, u_2, \ldots, u_6$, by identifying $v_5$ with $u_5$ (that is, $v_5=u_5$ in $G$).  Then $G$ is $5^*$-dense ($G$ is not $6$-dense since $v_1$ and $v_2$ have only $3$-common neighbors) and $G - u_6$ is $5^*$-dense as well.

\medskip

The next proposition gives us a sufficient condition for a graph $G$ on $n$ vertices to be $k$-dense.

\medskip

\bpr\quad If $G$ is a graph of order $n$ and $\delta(G) \geq \frac{n+k}{2}-1$, then $G$ is $k$-dense.\epr

\noindent
\pf\quad
Let $u$ and $v$ be any two adjacent vertices in $G$.  Since $\delta(G) \geq \frac{n+k}{2}-1$, there are at least $2\left[\frac{n+k}{2}-1\right]-2 = n+k-4$ edges between $u$ and $v$ and the remaining $n-2$ vertices in the graph.  Let $x$ be the number of common neighbors of $u$ and $v$.  Then $n+k-4-x \le n-2.$ This implies that $x \ge k-2$, and so $u$ and $v$ have at least $k-2$ common neighbors.\hfill\qed

\medskip

We now look at the edge connectivity of a $k$-dense connected graph.

\medskip

\bpr\label{PROP:edgeconnectivity}\quad  If a graph $G$ is $k$-dense and connected, then $G$ is at least $(k-1)$-edge-connected.\epr

\noindent
\pf\quad Suppose that $G$ is $k$-dense and connected.  If $k=2$, then since $G$ is connected, $G$ is $1$-edge-connected.  If $k \geq 3$, consider any set of $k-2$ edges $e_1, e_2, \ldots, e_{k-2}$.  We claim that $G-\{e_1, e_2, \ldots, e_{k-2}\}$ is connected.  Let $u,v \in V(G)$.  Since $G$ is connected, there exists a $u-v$ path in $G$, say $u=u_0, u_1, u_2, \ldots, u_p=v$.  For each $i$, $0 \leq i \leq p-1$, if $u_iu_{i+1} = e_j$ for some $j$, then notice that $u_i$ and $u_{i+1}$ have at least $k-2$ common neighbors. Since there are only $k-3$ remaining edges $e_j$, one of the common neighbors of $u_i$ and $u_{i+1}$ in $G$ is still a common neighbor in $G-\{e_1, e_2, \ldots, e_{k-2}\}$, so we can replace the edge $u_iu_{i+1}$ with a path of length $2$ through this common neighbor. Thus, there is still a $u-v$ walk in $G-\{e_1, e_2, \ldots, e_{k-2}\}$.
\hfill\qed

\medskip

The converse of Proposition~\ref{PROP:edgeconnectivity} is not true. Moreover, there is no similar result for vertex connectivity. These are shown in the two results that follow.

\medskip

\bpr\quad  For every positive integer $a$, there is a graph that is $2a$-connected and $2a$-edge connected that is $2^*$-dense. \epr

\noindent
\pf\quad The cartesian product of $a$ copies of $C_4$ is a graph that is connected and contains no triangles. Hence, the graph is 2-dense but not 3-dense.  A minimum vertex cut consists of the $2a$ neighbors of a vertex, and a minimum edge cut consists of the $2a$ edges incident with a vertex.\hfill\qed

\bpr\quad  For every positive integer $k \geq 3$, there exists a connected graph $G$ that is $k$-dense but not $2$-connected. \epr

\noindent
\pf\quad Take two copies of $K_k$ and identify a vertex.  The resulting graph is $k$-dense yet has a cut-vertex.~\hfill\qed


\section{Bounds, realization, and characterizations} \label{section:general}

For any integer $k$ with $k\ge 2$, the complete graph $K_k$ is $k$-dense but not $(k+1)$-dense. Hence, we know that a $k^*$-dense graph exists for each integer $k\ge 2.$  A more interesting question is that of determining the values of $k$ and $n$ for which there is a graph on $n$ vertices that is $k^*$-dense.  We will answer this question in this section and then classify which graphs on $n$ vertices are $2^*$-dense, $(n-1)^*$-dense, and $n^*$-dense.

\medskip

We begin by giving bounds for $k$ in terms of the number of vertices of the graph. Since the edge multiplicity of every edge in a graph on $n$ vertices  is between $0$ and $n-2$, we have the following.

\bobs\label{OBS:kdensebounds} If $G$ is $k^*$-dense graph on $n$ vertices, then $2\le k\le n-2.$
\eobs

The next result shows that the bounds given in Observation~\ref{OBS:kdensebounds} are sharp and characterizes which graphs satisfy the given bounds.


\bpr\label{PROP:ndense}\quad Let $n$ be an integer such that $n\ge 3$ and let $G$ be a graph on $n$ vertices. The following are true.
\begin{enumerate}
\item[{\rm (a)}] The graph $G$ is $2^*$-dense if and only if $G$ has no isolated vertices and there is an edge $uv$ in $G$ such that $N(u)\cap N(v)=\emptyset.$
\item[{\rm (b)}] The graph $G$ is $n^*$-dense if and only if $G=K_n.$
\end{enumerate}
\epr

\noindent
\pf\quad The first part follows readily from the definition of  $k^*$-dense given in Definition~\ref{DEFN:kdense} with $k=2.$

\medskip

Let us now establish the second part. Since any two vertices in $K_n$ are adjacent, it follows that every edge in $K_n$ has edge multiplicity $n-2$ and so $K_n$ is $n$-dense. By definition, no graph on $n$ vertices is $(n+1)$-dense. It follows that $K_n$ is $n^*$-dense. We now establish the converse by proving its contrapositive. Let $G$ be a graph on $n$ vertices such that $G\ne K_n.$ Thus, there exists two vertices $u$ and $v$ that are not adjacent in $G.$
Since any graph on $n$ vertices with isolated vertices is not $n$-dense, we may assume that $u$ is adjacent to some vertex $w$ in $G.$ Since $v\notin N(u)$, it follows that the edge $uw$ has edge multiplicity at most $n-3.$ Thus, $G$ is not $n$-dense and hence, is not $n^*$-dense.

\hfill\qed

\medskip

We are now ready to determine for which values of $k$ and $n$ there exists a graph $G$ on $n$ vertices that is $k^*$-dense.

\bpr\label{PROP:bounds-dense}\quad If $k$ and $n$ are integers such that $2\le k\le n$, then there is a $k^*$-dense graph on $n$ vertices.
\epr

\noindent
\pf\quad  Let $q$ and $r$ be the unique pair of integers such that $n=kq+r$ where $0\le r\le k-1.$ If $r=0$, let $G=qK_k$ and if $0<r\le k-1$, let $G=(q-1)K_k\cup (K_k+K_r).$  One can verify that in both cases, $G$ is a $k^*$-dense graph on $n$ vertices. \hfill\qed

\medskip

The next result characterizes the graphs on $n$ vertices that are $(n-1)^*$-dense.

\bpr\label{PROP:(n-1)dense}\quad If $n$ is an integer such that $n\ge 3$, then $G$ is $(n-1)^*$-dense if and only if $G \cong K_n-e.$
\epr

\noindent
\pf\quad Observe that every edge in $K_n-e$ has edge-multiplicity at least $n-3.$ Thus, $K_n-e$ is $(n-1)$-dense but not $n$-dense and so $K_n-e$ is $(n-1)^*$-dense. From Proposition~\ref{PROP:ndense}, $K_n$ is $n^*$-dense. Now, let $G$ be a graph on $n$ vertices such that $G\ne \{K_n,K_n-e\}.$ It follows that there exist distinct $uv, wx\in E(\overline{G})$ for some vertices (not necessarily all distinct) $u,v,w,x\in V(G).$ We may assume that $G$ has no isolated vertices for if $G$ does, then $G$ is not $(n-1)$-dense. Suppose first that $\{u,v\}\cap\{w,x\}\ne\emptyset$, say $u=w.$ It follows that $v\ne x$ and each of $v$ and $x$ is not in $N(u).$ Moreover, every edge incident to $u$ will have edge multiplicity at most $n-4.$ Thus, $G$ is not $(n-1)$-dense. Suppose now that $u, v, w$ and $x$ are distinct vertices; that is, $\{u,v\}\cap\{w,x\}=\emptyset.$ If $w\in N(u)$, then the edge $uw$ will have edge multiplicity at most $n-4$ and $G$ is not $(n-1)$-dense. Using similar arguments, one can show that $G$ is not $(n-1)$-dense if either $x\in N(u)$ or $\{w,x\}\cap N(v)\ne\emptyset.$ Thus, we may assume that $G[\{u,v,w,x\}]=\overline{K_4}.$ Since $G$ has no isolated vertices, there is an edge incident to $u$ in $G.$ Now every edge incident to $u$ in $G$ has edge multiplicity at most $n-5$ and so $G$ is not $(n-1)$-dense. Consequently, $K_n-e$ is the only graph on $n$ vertices that is $(n-1)^*$-dense.\hfill\qed


\section{Minimum number of edges in $2^*$-dense, $3^*$-dense and $4^*$-dense connected graphs on $n$ vertices}
\label{section:k=3,k=4}

In  Section~\ref{section:general}, we characterized all graphs on $n$ vertices that are $k^*$-dense for $k=2$,   $k=n-1$ or $k=n$.  We now investigate lower bounds for the number of edges in a graph on $n$ vertices that are $k^*$-dense for some small values of $k.$  We begin with a definition.


\begin{defn}\label{DEFN:kdensemin}\quad For $k,n\in\mathds{N}$, let $$e(k,n)=\min\{|E(G)| : G\ {\rm is\ a\ connected\ graph\ on\ } n {\rm\  vertices \ that\ is}\ k^*{\rm -dense}\}.$$
\end{defn}

Let us now consider a connected graph $G$ on $n$ vertices. Proposition~\ref{PROP:ndense} says that $G$ is $2^*$-dense if and only if $G$ has an edge that does not belong to a triangle in $G.$ Observe that any connected graph is $2$-dense and that trees have the minimum number of edges among all connected graphs on $n$ vertices. Since trees do not have triangles, we have the following the result.



\bth\label{Thm:2densemin}\quad For all $n\in\mathds{N}$ with $n\ge 2$, the minimum number of edges in a connected $2^*$-dense graph is  $e(2,n)=n-1,$ and this is achieved by any tree on $n$-vertices.
\eth


The next result tells us what $e(3,n)$ is for all $n\ge 3.$
\medskip

\bth\label{Thm:3densemin}\quad For all $n\in\mathds{N}$ with $n\ge 3$, the minimum number of edges in a connected $3^*$-dense graph is  $e(3,n)=\left\lceil\frac{3}{2}(n-1)\right\rceil.$
\eth

\noindent
\pf\quad   We consider two cases according to whether $n$ is odd or even.

\medskip

\noindent
\textit{Case $1:$ $n$ is odd.}

\medskip

Let $n=2\ell+1$ where $\ell\in\mathds{N}$. We need to show that $e(3, 2\ell+1)=3\ell.$ Consider $\ell$ copies of $K_3$ all identified at one vertex $w$.  Since this is a $3^*$-dense graph with $n$ vertices and $3\ell$ edges, we know $e(3, 2\ell+1)\geq 3\ell$.

\medskip

By Proposition~\ref{PROP:ndense}, the result is true when $\ell=1.$ Let us assume that $e(3, 2r+1)=3r$ for some $r\in\mathds{N}$ with $r\ge 1.$ We need to show that $e(3,2r+3)=3(r+1)=3r+3.$ Assume to the contrary that $e(3, 2r+3)\le 3r+2.$ Let $G$ be a connected graph with $2r+3$ vertices that is $3^*$-dense such that $G$ has $e(3, 2r+3)$ edges.

\medskip

We first show that $G$ has at least two vertices of degree 2. Recall that by Proposition~\ref{Prop:2.2}, every vertex in $G$ has degree at least 2. If $G$ has at most one vertex of degree 2, then $e(3, 2r+3)=|E(G)|\ge \frac{1}{2}(2+3(2r+2))=3r+4.$ But this contradicts our assumption that $e(3, 2r+3)\le 3r+2.$

\medskip

Now let $x,y\in V(G)$ such that $\deg_G x=2=\deg_G y.$ Let us first assume that $xy\in E(G).$ It follows from Definition~\ref{multiplicity_kdense} that $|N(x)\cap N(y)|=1$, say $N(x)\cap N(y)=\{u\}.$ Since $\deg_G x=2=\deg_G y$, we know that each of the edges $xy, xu$, and $yu$ has edge multiplicity 1 in $G$ (that is, $x$ and $y$ belong to only one triangle in $G$, namely $G[x,y,u]$).This means that $G-\{x,y\}$ is a connected $3^*$-dense graph with $2r+1$ vertices and $e(3, 2r+3)-3$ edges. From our assumption that $e(3, 2r+3)\le 3r+2$ it follows that $G-\{x,y\}$ has at most $(3r+2)-3=3r-1$ edges. But this means that $e(3, 2r+1)\le 3r-1$ which contradicts our inductive assumption that $e(3, 2r+1)=3r.$

\medskip

Let us now assume there exists no vertices $x$ and $y$ which are adjacent in $G$ if $\deg_G x=2=\deg_G y.$ Observe that $G-\{x,y\}$ is a connected graph with $2r+1$ vertices that has at most $(3r+2)-4=3r-2$ edges. Since$G-\{x,y\}$ has less than $e(3, 2r+1)$ edges, it follows that $G-\{x,y\}$ is not $3^*$-dense.

\medskip

Since $G$ is $3^*$-dense, each of $x$ and $y$ must be in exactly one triangle, say $x$, $u$ and $v$ form a triangle and $y$, $w$, and $z$ form a triangle.  Furthermore, since $G-\{x,y\}$ is not $3^*$-dense, either $uv$ or $wz$ has edge multiplicity $0$ in $G-\{x,y\}$, say without loss of generality that $uv$ has edge multiplicity $0$ in $G-\{x,y\}$.  Let $H$ be the graph obtained from $G-x$ by deleting the edge $uv$ and then identifying the vertices $u$ and $v$.  Since $uv$ had edgemultiplicity $0$, this will not affect the multiplicity of any other edge in $G-x$. Note that, by construction, $H$ is a $3^*$-dense graph on $2r+1$ vertices that has $e(3,2r+3)-3$ edges.  This means that $H$ has at most $(3r+2)-3=3r-1$ edges and so $e(3,2r+1) \leq 3r-1$.  But again, this contradicts our inductive assumption that $e(3,2r+1)=3r$.  We therefore conclude that $e(3,2r+3) = 3r+3$.  Consequently, $e(3,2\ell+1)=3\ell$ for all $\ell \in \mathds{N}.$

\medskip

\noindent
\textit{Case $2:$ $n$ is even.}

\medskip

Let $n=2\ell$ where $\ell\in\mathds{N}$ with $\ell\ge 2.$ We need to show that $e(3, 2\ell)=3\ell-1.$  Consider $\ell-2$ copies of $K_3$ and one copy of $K_4 - e$ all identified at one vertex $w$.  This is a $3^*$-dense graph with $n$ vertices and $3\ell-1$ edges, so we have $e(3, 2\ell)\geq 3\ell-1$.

\medskip

 By Proposition~\ref{PROP:(n-1)dense}, the result is true when $\ell=2.$ Let us assume that $e(3, 2r)=3r-1$ for some $r\in\mathds{N}$ with $r\ge 2.$ We need to show that $e(3, 2r+2)=3(r+1)-1=3r+2.$ Assume to the contrary that $e(3, 2r+2)\le 3r+1.$ Let $G$ be a connected graph with $2r+2$ vertices that is $3^*$-dense such that $G$ has $e(3, 2r+2)$ edges. We first show that $G$ has at least one vertex of degree 2. If $\deg v\ge 3$ for all $v\in V(G)$, then $e(3, 2r+2)=|E(G)|=\frac{1}{2}(3(2r+2))\ge 3r+3.$ But this contradicts our assumption that $e(3, 2r+2)\le 3r+1.$ Let $x\in V(G)$ such that $\deg x=2.$ Consider the graph $G-x.$ Observe that $G-x$ is a connected graph with $2r+1$ vertices and $e(3,2r+2)-2$ edges. Thus, $G-x$ has at most $(3r+1)-2=3r-1$ edges. Since $e(3, 2r+1)=3r$ from Case 1 above, 
  we know that $G-x$ is not $3^*$-dense. This means that there exist vertices $w$ and $z$ such that $G[w,x,z]$ is a triangle in $G$ such that the edge $wz$ has edge multiplicity 0 in $G-\{x\}$ (that is, $N_{G-\{x\}}(w)\cap N_{G-\{x\}}(z)=\emptyset$). Let $H$ be the graph obtained from $G-\{x\}$ by deleting the edge $wz$ and then identifying the vertices $w$ and $z.$ Note that, by construction, $H$ is a $3^*$-dense graph on $2r$ vertices that has $e(3, 2r+2)-3$ edges. This means that $H$ has at most $(3r+1)-3=3r-2$ edges and so $e(3, 2r+1)\le 3r-2.$ But this contradicts our inductive assumption that $e(3,2r)=3r-1.$ We therefore conclude that $e(3, 2r+2)=3r+2.$ Consequently, $e(3, 2\ell)=3\ell-1$ for all $\ell\in\mathds{N}$ with $\ell\ge 2.$ \hfill\qed

\medskip


We now determine $e(4,n)$ for all $n\ge 4.$

\begin{thm}  \label{Thm:4densemin} For every positive integer $n \geq 4$, the minimum number of edges in a connected $4^*$-dense graph is
$$e(4,n) = \left\{\begin{array}{cl}
2n-2  & {\rm if}\  n \equiv 1\mod 3 \\
2n-1 & {\rm  otherwise.}\\
\end{array}
\right.$$
\end{thm}

\noindent
\pf\quad  For convenience in the proof, let
$$F(n) = \left\{\begin{array}{cl}
2n-2  & {\rm if}\  n \equiv 1\mod 3 \\
2n-1 & {\rm  otherwise.}\\
\end{array}
\right.,$$
 so that we wish to prove that $e(4,n)=F(n)$.  Notice that $F(n-3)=F(n)-6$ for all $n$.

First, we show that $e(4,n) \leq F(n)$ using an inductive construction.  For $n=4$, let $G=K_4$.  For $n=5$, let $G=K_5-e$.  For $n=6$, we have the graph $G = K_2+2K_2$ that is a $4^*$-dense graph with 6 vertices and 11 edges.  Suppose we have a graph $G$ with $n$ vertices and $F(n)$ edges.  Associate one vertex of the complete graph $K_4$ with a vertex of $G$ to obtain a $4^*$-dense graph with $n+3$ vertices and $F(n)+6=F(n+3)$ edges.

\medskip

Next, we must show that $e(4,n) \geq F(n)$.  We proceed by induction on $n$. For $n = 4$ and $n=5$, the result follows from Propositions~\ref{PROP:ndense} and~\ref{PROP:(n-1)dense}.  Suppose we have a $4^*$-dense graph on $6$ vertices.  From Proposition~\ref{Prop:2.2}, we know that every vertex has degree at least $3$, and furthermore, any vertex with degree 3 is in a subgraph isomorphic to $K_4$.  If every vertex has degree $4$ or more, then the graph has $12$ edges.  Suppose there is a vertex $u$ with $\deg u = 3$.  Then $u$ lies in a subgraph isomorphic to $K_4$, so there are vertices $v$, $w$ and $x$ such that $u$, $v$, $w$ and $x$ induce a complete subgraph, with $6$ edges.  The remaining two vertices in the graph, say $y$ and $z$, each have degree at least $3$.  If $yz$ is an edge, then there must be at least $4$ more edges between $\{u,v,w,x\}$ and $\{y,z\}$, for a total of $6+1+4=11$ edges.  If $yz$ is not an edge, then there are at least $6$ edges between $\{u,v,w,x\}$ and $\{y,z\}$, for a total of at least $6+6=12$ edges.  Thus, a connected $4^*$-dense graph on $6$ vertices must have at least $11$ edges.

\medskip

Consider an integer $n \geq 7$.  Assume, for all $k$ such that $4 \leq k < n$, we have $e(4,k)=F(k)$.  We wish to show that $e(4,n)=F(n)$.  Suppose, to the contrary, that there is a connected $4^*$-dense graph $G$ with $n$ vertices and at most $F(n)-1$ edges.  Thus, $G$ has at most $2n-2$ edges.  From Proposition~\ref{Prop:2.2}, the minimum degree in $G$ is at least $3$.  If every vertex in $G$ has degree at least $4$, then $G$ would have at least $\frac{4n}{2}=2n$ edges, which is a contradiction.  Thus, there is some vertex $u$ in $G$ such that $\deg u = 3$. Let $N(u) = \{v,w,x\}$. From Proposition~\ref{Prop:2.2}, we know that $\{u, v, w, x\}$ induces a complete subgraph.

\medskip

Notice that in $G-u$, each of the edges $vw$, $wx$, and $xv$ is in the triangle $v, w, x$, and hence has multiplicity at least 1.  If all three edges have edge multiplicities equal to 2, then $G-u$ is $4^*$-dense.

\medskip

\noindent
\textit{Case $1:$ $G-u$ is $4^*$-dense.}

\medskip

In this case, $G-u$ has $n-1$ vertices and at most $(F(n)-1)-3=F(n)-4$ edges.  If $n \equiv 0 \mod 3$, then $G-u$ has at most $(2n-1)-4=2n-5$ edges.  Since $n-1 \equiv 2 \mod 3$, we have $F(n-1)=2(n-1)-1=2n-3$.  If $n \equiv 1 \mod 3$, then $G-u$ has at most $(2n-2)-4 =2n-6$ edges.  Since $n-1 \equiv 0 \mod 3$, we have $F(n-1) = 2(n-1)-1=2n-3$.  If $n \equiv 2 \mod 3$, then $G-u$ has at most $(2n-1)-4=2n-5$ edges.  Since $n-1 \equiv 1 \mod 3$, we have $F(n-1) = 2(n-1) - 2 = 2n-4$.  In each case, $G-u$ has fewer than $F(n-1)$ edges, which contradicts the inductive hypothesis that $F(k)=e(4,k)$ for all $k$ with $4 \leq k < n$.

\medskip

\noindent
\textit{Case $2:$ The edges $vw$, $vx$, and $wx$ all have edge multiplicity 1 in $G-u$.}

\medskip

Notice that the edge multiplicity 1 in each case comes from the triangle $v, w, x$.  Thus, we can delete the edges $vw$, $vx$, and $wx$ without changing the edge multiplicity of any of the other edges in $G-u$ and then identify all three vertices to obtain a new graph $H$ with $n-3$ vertices and at most $(F(n)-1)-6=F(n)-7$ edges.  Notice that $H$ is connected and $4^*$-dense by construction.  Since $F(n-3) = F(n)- 6$ for all $n \geq 7$, the graph $H$ has order $n-3$ and fewer than $F(n-3)$ edges, which contradicts the inductive hypothesis.

\medskip

\noindent
\textit{Case $3:$  Exactly one of the edges $vw$, $vx$, and $wx$, say $vx$, has edge multiplicity 1 in $G-u$.}

\medskip

We may assume that $vw$ and $wx$ each have edge multiplicity at least $2$ in $G-u$.  Thus, $v$ and $w$ have at least one common neighbor, say $y$, other than $x$, and $w$ and $x$ have at least one common neighbor, say $z$, other than $v$.  Notice that if $y=z$, then $y$ is a common neighbor of $v$ and $x$, which contradicts our assumption that $vx$ has edge multiplicity 1 (that $w$ is the only common neighbor of $v$ and $x$).

\medskip

We may delete the edge $vx$ from $G-u$ without changing the edge multiplicity of any other edge except for $vw$ and $wx$, and then identify the vertices $v$ and $x$, merging the two edges $vw$ and $wx$ into one new edge.  The endpoints of this new edge have at least two common neighbors, $y$ and $z$.  Thus, the resulting graph $H$ is $4$-dense and connected, with $n-2$ vertices and at most $(F(n)-1)-5=F(n)-6$ edges.

If $n \equiv 0 \mod 3$, then $F(n)-6 = (2n-1)-6 = 2n-7$.  Since $n-2 \equiv 1 \mod 3$, we have $F(n-2)=2(n-2)-2=2n-6$.
If $n \equiv 1 \mod 3$, then $F(n)-6 = (2n-2)-6 = 2n-8$.  Since $n-2 \equiv 2 \mod 3$, we have $F(n-2)=2(n-2)-1=2n-5$.
If $n \equiv 2 \mod 3$, then $F(n)-6 = (2n-1)-6 = 2n-7$.  Since $n-2 \equiv 0 \mod 3$, we have $F(n-2)=2(n-2)-1=2n-5$.
In each case, we have a connected $4^*$-dense graph with $n-2$ vertices and fewer than $F(n-2)$ edges, which contradicts the inductive hypothesis.

\medskip

\noindent
\textit{Case $4:$  Exactly two of the edges $vw$, $vx$, and $wx$, say $vx$ and $vw$, have edge multiplicity 1 in $G-u$.}

\medskip

We may delete $vw$ and $vx$ in $G-u$ without changing the edge multiplicity of any other edge in the graph except for $wx$, and then identify $v$ with $w$, to ensure that the resulting graph $H$ is connected.  The edge multiplicity of $wx$ in $H$ is at least 1.  If the edge multiplicity of $wx$ is at least 2, then $H$ is a connected $4^*$-dense graph with $n-2$ vertices and at most $(F(n)-1)-5=F(n)-6$ edges.  This contradicts the inductive hypothesis, as in Case $2$.

\medskip

We may assume that $wx$ has edge multiplicity 1 in $H$.  Let $y$ be the unique common neighbor of $w$ and $x$ in $H$.  Since $yw$ and $yx$ each have multiplicity at least $2$ in $G$ and in $H$, we know that $x$ and $y$ have a common neighbor $z$ with $z \neq w$, and $w$ and $y$ have a common neighbor $z'$ such that $z' \neq x$.  If $z=z'$, then $w$ and $x$ have common neighbors $y$ and $z$, which contradicts our assumption that $wx$ has multiplicity 1.  If $z \neq z'$, then we will delete the edge $wx$ and identify the vertices $w$ and $x$, merging the edges $xy$ and $wy$ into a single new edge.  The endpoints of this new edge have common neighbors $z$ and $z'$, so the new edge has multiplicity at least 2.  We also add the edge $zz'$.  Notice that $z$ and $z'$ have at least two common neighbors, $y$ and the vertex formed from identifying $w$ and $x$, so $zz'$ has edge multiplicity at least 2.  The edge multiplicity of the remaining edges is not decreased, so the new graph $J$ has $n-3$ vertices and at most $(F(n)-1)-7+1=F(n)-7$ edges.  We can readily check that $F(n-3)=F(n)-6$ for all $n$, so this contradicts our inductive hypothesis.
\hfill\qed

\medskip

Notice that in Theorem~\ref{Thm:2densemin}, Theorem~\ref{Thm:3densemin}, and Theorem~\ref{Thm:4densemin}, $$e(k,n) = q{k\choose 2}+{r\choose 2}+r(k-r),$$ where $q$ and $r$ are the unique integers such that $n-1=q(k-1)+r$ and $0\le r<k-1.$ The next section shows that this need not be the case for $k \ge 8.$ We thus have the following conjecture.

\bcon\quad Let $G$ be a $k$-dense connected graph with $n$ vertices, where $k\le 7.$ If $q$ and $r$ are the unique integers such that $(n-1) = q(k-1)+r$ and $0 \leq r < k-1$,  then $$e(k,n) = q{k\choose 2}+{r\choose 2}+r(k-r).$$
\econ

\section{Upper bounds on the minimum number of edges in $k^*$-dense graphs on $n$ vertices} \label{section:k<7}

In Section~\ref{section:general}, we characterized all graphs on $n$ vertices that are $2^*$-dense, $(n-1)^*$-dense, or $n^*$-dense, while in Section~\ref{section:k=3,k=4},  we found the minimum number of edges in connected graphs on $n$ vertices that $2^*$-dense, $3^*$-dense or $4^*$-dense.   In this section, we present a result that gives us an upper bound on the minimum number of edges in graphs (not necessarily connected) on $n$ vertices that are $k^*$-dense, for all $k$ such that $2\le k\le n.$


\bpr\label{PROP:kdensemingeneral}\quad Let $k$ and $n$ be integers such that $2\le k\le n$. If $q$ and $r$ are the unique integers such that $n=kq+r$ and $0\le r<k$, then there exists a $k^*$-dense graph on $n$ vertices that has
$$ q \binom{k}{2} + \binom{r}{2} + r (k-r)$$
edges.
\epr	

\noindent
\pf\quad The graph $G$ ($G=qK_k$ if $r=0$ and $G=(q-1)K_k\cup K_k+K_r$ if $0<r<k$) given in the proof of Proposition~\ref{PROP:bounds-dense} is a $k^*$-dense graph on $n$ vertices that has
$q \binom{k}{2} + \binom{r}{2} + r (k-r)$ edges.\hfill\qed

\medskip

In Propostion~\ref{PROP:kdensemingeneral}, we did not require the $k^*$-dense graph to be connected. If we consider connected graphs only, then we have the following.

\bpr\label{PROP:kdenseminconnected}\quad Let $k$ and $n$ be integers such that $2\le k\le n$. If $q$ and $r$ are the unique integers such that $n-1=(k-1)q+r$ and $0\le r<k-1$, then there exists a $k^*$-dense connected graph on $n$ vertices that has
$$ q \binom{k}{2} + \binom{r}{2} +r (k-r)$$
edges.
\epr

\noindent
\pf\quad Let $H_1=K_k+K_r$ and for $i=2, 3, \cdots, q$, let $H_i=K_k.$ For $i=1, 2, \cdots, q$, let $v_i$ be a vertex of $H_i.$ If $G$ is the graph obtained by identifying the vertices $v_1, v_2,\cdots, v_q$, then $G$ is a $k^*$-dense connected graph on $n$ vertices that has $q\binom{k}{2}+\binom{r}{2}+r(k-r)$ edges. \hfill\qed

\medskip

It follows from Proposition~\ref{PROP:kdenseminconnected} that, in general, $e(k,n)\le q \binom{k}{2} + \binom{r}{2} + r (k-r)$, where $k, n, q,$ and $r$ are as given above. However, for certain values of $k$ and $n$, we can improve this upper bound for $e(k,n)$ as we now show.

\begin{obs}\label{OBS:alternative}
Let $$H=K_{n-k}\cup \left(\frac{k}{2}\right)K_2$$ if $k$ is even and let $$H=K_{n-k}\cup\left( \frac{k-1}{2}\right)K_2\cup K_1$$ if $k$ is odd, where $n-k \geq 2$. The graph $G=\overline{H}$ is a connected graph on $n$ vertices that is $k^*$-dense.
\end{obs}

Let us consider particular examples of the graphs described in Observation~\ref{OBS:alternative}. If $n=26$ and $k=23$, then $H=K_3\cup 11 K_2\cup K_1$; that is $H$ is the union of one triangle, 11 independent edges, and an isolated vertex and the graph $G$ is the complement of $H.$ In this construction, $G$ is a connected graph with $26$ vertices that is $23^*$-dense and has $${26\choose 2}-3-11=325-13=311$$ edges. Note that $25=1\cdot 22+3.$ Hence, the construction of Proposition~\ref{PROP:kdenseminconnected} gives us a  connected graph with $26$ vertices that is $23^*$-dense and has $$1\cdot {23\choose 2}+{3\choose 2 }+3\cdot(23-3)=316$$ edges. As another example, if $n=26$ and $k=24$, then $H=13 K_2$; that is $H$ is the union of 13 independent edges and the graph $G$ is the complement of $H.$ In this construction, $G$ is a connected graph with $26$ vertices that is $24^*$-dense and has $${26\choose 2}-13=325-13=312$$ edges. Note that $25=1\cdot 23+2.$ Hence, the construction of Proposition~\ref{PROP:kdenseminconnected} gives us a  connected graph with $26$ that is $24^*$-dense and has $$1\cdot {24\choose 2}+{2\choose 2}+2\cdot(24-2)=321$$ edges.\\


For many values of $k$ and $n$, the upper bound for $e(k,n)$ given by Proposition~\ref{PROP:kdenseminconnected} is better than the one given by Observation~\ref{OBS:alternative} (say $n=10, k=6$ or $n=10,k=7$). But the examples above  show that the opposite is true for some values of $k$ and $n$ ($n=26,k=23$ or $n=26,k=24$).


\section{Maximum number of edges in $k^*$-dense connected graphs on $n$ vertices}\label{section:upper_bound}


In this section, we determine the maximum number of edges in a $k^*$-dense connected graph on $n$ vertices. We start by giving a definition that is analogous to $e(k,n).$

\begin{defn}\label{DEFN:kdensemax}\quad For $k,n\in\mathds{N}$, let $$E(k,n)=\max\{|E(G)| : G\ {\rm is\ a\ connected\ graph\ on\ } n\ {\rm vertices\ that\ is\ } k^*{\rm -dense}\}.$$
\end{defn}

\medskip



\begin{thm} \label{Thm:upper}
 If $k$ and $n$ are integers such that $2\le k\le n$, then $E(k, n) = n+k-3+ \left(\begin{array}{cc} n-2 \\ 2 \end{array} \right).$ \end{thm}

\noindent
\pf\quad Suppose first that $n \geq k + 2$ and $k \geq 2$ and let $G$ be a graph on $n$ vertices with $n+k-2+ \left(\begin{array}{cc} n-2 \\ 2 \end{array} \right)$ edges.  Let $u$ and $v$ be any two adjacent vertices in $G$.  Note that there are at most $\left(\begin{array}{cc} n-2 \\ 2 \end{array} \right)$ edges that are not incident with either $u$ or $v$ (consider the edges in the subgraph induced by the vertex set $V(G)-\{u,v\}.$)  Thus, there are at least $n+k-3$ edges between $u$ and $v$ and the vertex set $V(G)-\{u,v\}.$  By the pigeonhole principle, $u$ and $v$ must have at least $k-1$ common neighbors.  Since this is true for any adjacent pair $u$ and $v$, it follows that $G$ is $(k+1)$-dense. This means that $G$ is not $k^*$-dense. It follows that $E(k,n)\le n+k-3+ \displaystyle{n-2\choose 2}.$

\medskip

To see that this is sharp, consider the graph $G$ formed by starting with a complete graph $K_{n-2}$.  Now add two vertices $u$ and $v$, with an edge between $u$ and $v$.  Partition the vertices of the $K_{n-2}$ into three sets $A$, $B$, and $C$ such that $\displaystyle{|A|=\left\lfloor\frac{n-k}{2}\right\rfloor}$, $|B|=k-2$ and $\displaystyle{|C|=\left\lceil\frac{n-k}{2}\right\rceil}$.  Join $u$ to every vertex in $A$ and $B$ and join $v$ to every vertex in $B$ and $C$.  We can check that $G$ has exactly $n+k-3+ \left(\begin{array}{cc} n-2 \\ 2 \end{array} \right)$ edges.  Since $u$ and $v$ are adjacent and have exactly $k-2$ common neighbors, specifically the vertices of $B$, we know that $G$ is not $(k+1)$-dense.  Now, any two vertices distinct from $u$ and $v$ are in the subgraph $K_{n-2}$, so they have at least $n-4 \geq k-2$ common neighbors. The vertices $u$ and $a\in N_G(u)$, where $a \neq v$, have at least $\displaystyle{\left\lfloor\frac{n-k}{2}\right\rfloor+k-3 \geq \frac{2}{2}+k-3 = k-2}$ common neighbors (similarly for the vertices $v$ and $b$ where $b\in N_G(v)$ and $b \neq u$).  Hence, $G$ is $k$-dense. Consequently, $G$ is $k^*$-dense.

Now, if $n=k$, then $n+k-3+ \left(\begin{array}{cc} n-2 \\ 2 \end{array} \right) = \left(\begin{array}{cc} k\\ 2 \end{array} \right) $.  We know that the only $k^*$-dense graph on $k$ vertices is $K_k$ by Proposition~\ref{PROP:ndense}.  If $n=k+1$, then $n+k-3+ \left(\begin{array}{cc} n-2 \\ 2 \end{array} \right)=\left(\begin{array}{cc} k+1 \\ 2 \end{array} \right) -1$ and by Proposition~\ref{PROP:(n-1)dense}, the only $k^*$-dense graph on $k+1$ vertices is $K_{k+1}-e$.  \hfill\qed


\section{Complete Results for $2^*$-dense,  $3^*$-dense and $4^*$-dense Connected Graphs}\label{Section:complete k=3, 4}

By Theorem~\ref{Thm:2densemin}  and Theorem~\ref{Thm:upper}, we have the bounds in the following corollary.

\bco\label{COR:2dense}\quad
Let $n$ be an integer with $n\ge 4.$ If $G$ is $2^*$-dense graph with $n$ vertices and $m$ edges, then
	\begin{equation} \label{eqn:2dense}
	n-1\le m \le  {n-2\choose 2}+n-1.
	\end{equation}
	Moreover,  for every integer $a$  with $n-1 \le a\le {n-2 \choose 2}+n-1,$ there is  a $2^*$-dense graph $G$ with $n$ vertices and $a$ edges.
\eco

\noindent
\pf \quad Equation~\ref{eqn:2dense} follows from  Theorem~\ref{Thm:2densemin}  and Theorem~\ref{Thm:upper}.

\medskip

For the realization result, let $u_1$ and $u_2$ be two adjacent vertices of a tree on $n$ vertices, such that $\deg u_1 = 1$, then add $a-n+1$ more edges not incident to $u_1$.  Since $m(u_1, u_2) = 0$, this graph is $2$-dense. \hfill\qed

\medskip

By Theorem~\ref{Thm:3densemin}  and Theorem~\ref{Thm:upper}, we have the bounds in the following corollary.

\bco\label{COR:3dense}\quad
Let $n$ be an integer with $n\ge 4.$ If $G$ is $3^*$-dense graph with $n$ vertices and $m$ edges, then
	\begin{equation} \label{eqn:3dense}
	\frac{3}{2}(n-1)\le m \le  {n-2\choose 2}+n.
	\end{equation}
	Moreover,  $\forall a$  with $\frac{3}{2}(n-1)\le a\le {n-2 \choose 2}+n,$ there is  a $3^*$-dense graph $G$ with $n$ vertices and $a$ edges.
\eco

\noindent
\pf \quad Equation~\ref{eqn:3dense} follows from  Theorems~\ref{Thm:3densemin}  and~\ref{Thm:upper}.

\medskip

For the realization result, consider first an odd $n$.  Then the graph $G$ obtained from

\begin{enumerate}
	
	\item $\frac{n-1}{2}$ copies of $K_3$ all identified at one vertex $w$ (so that $\deg w = n-1$ and all other vertices have degree $2$),

	\item if $u_1$ and $u_2$ are two adjacent vertices of degree $2$, then add $a-\frac{3(n-1)}{2}$ more edges
such that $N(u_1) \cap N(u_2) = \{w\}$.


\end{enumerate}

Since $m(u_1, u_2) = 1$, it follows that  $G$ is an $n$-vertex graph that is $3^*$-dense, with $a$ edges.\\
For even $n$, use the same construction as above, but replace one of the triangles by a $K_4-e$, such that $\deg w = n-1$ is still valid. \hfill\qed

\medskip

Similarly, by Theorem~\ref{Thm:4densemin} and Theorem~\ref{Thm:upper}, we have the bounds in the following corollary.

 \bco\quad Let $n$ be an integer with $n\ge 4.$ If $G$ is $4^*$-dense graph with $n$ vertices and $m$ edges, then:

	\begin{equation} \label{eqn:4dense}
\left\{\begin{array}{cl}
2n-2  & {\rm if}\  n \equiv 1\mod 3 \\
2n-1 & {\rm  otherwise}\\
\end{array}
\right.
\le m
\le  n+1+ \left(\begin{array}{cc} n-2 \\ 2 \end{array} \right).
\end{equation}
	Moreover,  for all values $a$  between the lower and upper bounds of Equation~\ref{eqn:4dense}, there is  a $4^*$-dense graph $G$ with $n$ vertices and $a$ edges.
\eco
\pf \quad Equation~\ref{eqn:4dense} follows from  Theorems~\ref{Thm:4densemin}  and~\ref{Thm:upper}.

\medskip

For the realization result, consider first  $n \equiv 1 \mod 3$.  Then the graph $G$ is obtained from

\begin{enumerate}

	\item $\frac{n-1}{3}$ copies of $K_4$ all identified at one vertex $w$ (so that $\deg w = n-1$ and all other vertices have degree $3$),

	\item if $u_1, u_2$ and $u_3$ are three mutually adjacent vertices of degree $3$, then add $a-{ 4\choose 2}\frac{(n-1)}{3}$ more edges such that $N(u_1) \cap N(u_2)  = \{u_3, w\}$.
\end{enumerate}

Since $m(u_1, u_2) = 2$, it follows that  $G$ is an $n$-vertex graph that is $4^*$-dense, with $a$ edges. If $n \equiv 2 \mod 3$, use the same construction as above, but replace one of the copies of $K_4$ by a $K_5-e$ such that $\deg w = n-1$ is still valid. If $n \equiv 0 \mod 3$, use the same construction as above, but replace two of the copies of $K_4$ by two copies of $K_5-e$ such that $\deg w = n-1$ is still valid. ~\hfill\qed


\section{Conclusion}

In this paper, we introduced the topic of $k$-dense to the graph theory community, and we studied the minimum number of edges, $e(k,n)$, and the maximum number of edges, $E(k,n)$, that a connected $k^*$-dense ($k$-dense but not $(k+1)$-dense) graph of order $n$ can have. We established a formula for $E(k,n)$ and found an upper bound for $e(k,n)$ for all $2\le k \le n$ that is sharp for small values of $n.$  For small values of $k$ (namely $2, 3, 4$), we showed that there exists a $k^*$-dense connected graph $G$ on $n$ vertices that has $m$ edges for each $m$ satisfying $e(k,n)\le m\le E(k,n).$ Using a computer to assist us in our investigation (see Appendix), we believe that the formula for $e(k,n)$ that we found for $k=2,3,4$ generalizes to all $k\le 7.$  However, for $k \ge 8$ the formula for $e(k,n)$ that we found does not hold.  We conclude the paper with the following open problem.

\bprob\quad Let $k$ and $k$ be integers such that $5\le k\le n.$ Determine $e(k,n)$; that is, find the minimum number of edges that a graph $G$ on $n$ vertices can have if $G$ is $k^*$-dense. \eprob



\newpage

\section*{Appendix: Sage Code}
\label{appdx}

\begin{figure}[h]
\centering
\includegraphics{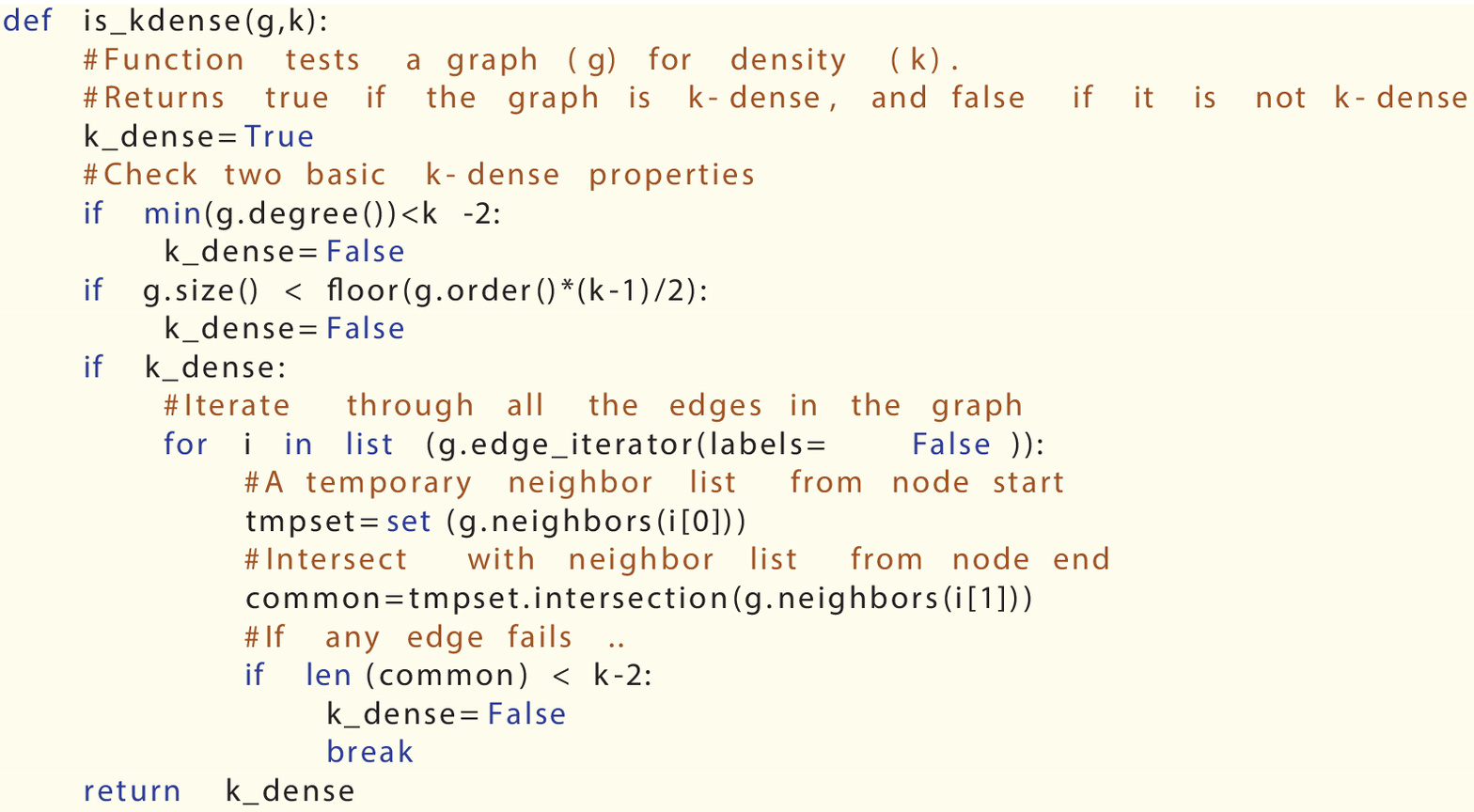}
\caption{Sage code to test a graph for k-density}
\label{fig:ktest}
\end{figure}

For checking graphs with $k \ge 5$, we used Sage\cite{sage}. Specifically, we generated graphs within the range of anticipated edge-counts for a given $n$, then tested each graph using the function found in Figure \ref{fig:ktest}.

\begin{figure}[h]
\centering
\includegraphics{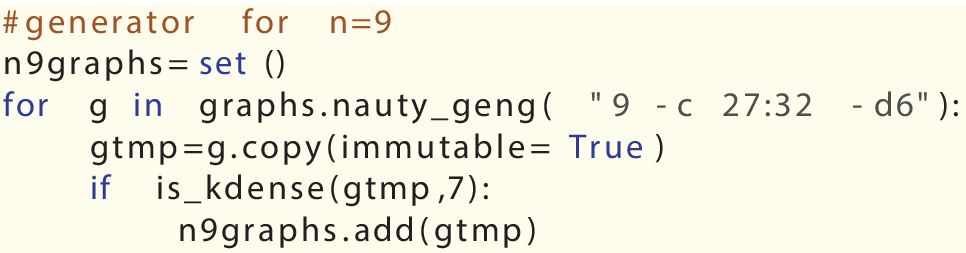}
\caption{Sample Sage code for testing graphs of $n=9$ for density of $k=7$.}
\label{fig:krun}
\end{figure}

\noindent The actual graphs we tested were generated using the `Nauty Geng' package created by McKay and Piperno \cite{nauty_geng}, which is included in Sage by permission. A sample of our code for testing $n=9$ is given in Figure \ref{fig:krun}.

\medskip

Note that while the sample code given here works, it has not been optimized or carefully coded for efficiency.

\end{document}